\documentclass[reqno, 11pt]{amsart}

\usepackage{amssymb}
\usepackage{latexsym}
\usepackage{amsmath}

\newcommand{\mcal}{\mathcal}

\makeatletter \@addtoreset{equation}{section}

\newtheorem{lemma}{Lemma}[section]
\newtheorem{theorem}[lemma]{Theorem}
\newtheorem{proposition}[lemma]{Proposition}
\newtheorem{corollary}[lemma]{Corollary}
\newtheorem{question}[lemma]{Question}
\theoremstyle{definition}
\newtheorem{definition}[lemma]{Definition}
\newtheorem{example}[lemma]{Example}
\newtheorem{remark}[lemma]{Remark}

\newcommand{\bdf}{\begin{definition}}
\newcommand{\edf}{\end{definition}}

\newcommand{\blem}{\begin{lemma}}
\newcommand{\elem}{\end{lemma}}

\newcommand{\bthm}{\begin{theorem}}
\newcommand{\ethm}{\end{theorem}}

\newcommand{\bpf}{\begin{proof}}
\newcommand\bpff[1]{\begin{proof}[#1]}
\newcommand{\epf}{\end{proof}}

\newcommand{\bprop}{\begin{proposition}}
\newcommand{\eprop}{\end{proposition}}

\newcommand{\bcor}{\begin{corollary}}
\newcommand{\ecor}{\end{corollary}}

\newcommand{\brem}{\begin{remark}}
\newcommand{\erem}{\end{remark}}

\newcommand{\bquest}{\begin{question}}
\newcommand{\equest}{\end{question}}

\newcommand{\bex}{\begin{example}}
\newcommand{\eex}{\end{example}}

\newcommand\0[1]{{\mathbb #1}}
\newcommand\eps{\varepsilon}
\newcommand\pG{\partial \Gamma}

\begin{document}
\title{The ratio set of the harmonic measure of a random walk on
a hyperbolic group}
\author{Masaki Izumi, Sergey Neshveyev and Rui Okayasu}
\address{}
\email{}
\date{}

\begin{abstract}
We consider the harmonic measure on the Gromov boundary of a
nonamenable hyperbolic group defined by a finite range random walk
on the group, and study the corresponding orbit equivalence
relation on the boundary. It is known to be always amenable and of
type~III. We determine its ratio set by showing that it is
generated by certain values of the Martin kernel. In particular,
we show that the equivalence relation is never of type III$_0$.
\end{abstract}

\maketitle

\section*{Introduction}
Given a discrete group $\Gamma$ and a probability measure $\mu$ on
$\Gamma$, one can construct a measure space $(X,\nu)$ with an
action of $\Gamma$ on it called the Poisson boundary. In many
cases this space can be identified with a geometric boundary of
$\Gamma$, and then the harmonic measure $\nu$ is the hitting
distribution of a random walk defined by $\mu$. E.g. if $\Gamma$
is a nonelementary hyperbolic group and $\mu$ is finitely supported
then the Poisson boundary coincides with the Gromov boundary~\cite{a}
(see also~\cite{k2}).
Conversely, a large class of measures on the Gromov boundary can be
obtained as harmonic measures for certain random walks
(though typically not of finite range)~\cite{cm1,cm2}.

It is known that the action of $\Gamma$ on $(X,\nu)$ is always
ergodic, amenable and, as soon as $(X,\nu)$ is nontrivial, of type
III. In particular, apart from the type III$_0$ case, the orbit
equivalence relation is completely determined by the ratio
set~\cite{cfw,kr} . There are a number of papers dealing with the
computation of the ratio sets for such actions. To name a few,
Sullivan~\cite{su0,su} proved that the recurrent part of an action
of a discrete conformal group on the sphere $S^d$ relative to the
Lebesgue measure is of type III$_1$. Spatzier~\cite{sp} showed
that if $\Gamma$ is the fundamental group of a compact connected
negatively curved manifold $M$ then the action of $\Gamma$ on the
sphere at infinity of the universal cover $\widetilde{M}$ of $M$
considered with the visibility measure class is again of type
III$_1$. For the case of free groups the ratio sets for certain
harmonic measures were computed by Ramagge and Robertson~\cite{rr}
and one of the authors~\cite{o}.

The result of Spatzier mentioned above is based on the work of
Bowen~\cite{b} who studied equivalence relations defined by Anosov
foliations. Bowen's computations are based in turn on the fact
that these equivalence relations are stably isomorphic to certain
equivalence relations arising from symbolic dynamical systems with
Gibbs measures. Since harmonic measures on the boundaries of free
groups are Gibbs measures~\cite{la,l,se}, Bowen's results
determine in principle the types of equivalence relations on the
boundaries of free groups. Moreover, in view of connections
between hyperbolic groups and symbolic dynamical systems (see
e.g.~\cite{cp}), one can hope that his results can be applied to
arbitrary hyperbolic groups. Indeed, following Bowen's strategy we
completely determine the ratio set for any nonelementary
hyperbolic group and the harmonic measure defined by a
nondegenerate finite range random walk on the group. In
particular, we show that the orbit equivalence relations we
consider are never of type III$_0$. Note that our result can be
thought of as a von Neumann algebra counterpart of the results of
Anantharaman-Delaroche~\cite{ad} and Laca and Spielberg~\cite{ls}
saying that the crossed product C$^*$-algebras arising from
boundary actions of hyperbolic groups are purely infinite.

We cannot directly apply Bowen's result on Gibbs measures to our
problem and we need to produce hyperbolic group versions of
various statements he used. One of the most important ingredients
of our proof is Ancona's result on almost multiplicativity of the
Green kernel, which was the key observation for identification of
the Martin boundary and the Gromov boundary ~\cite{a}. The
H\"older continuity of the Martin kernel also plays an essential
role and to establish it we follow Ledrappier's argument~\cite{l}
in the case of free groups. In the course of our proof, we also
show a vanishing theorem for certain first cohomology of the
boundary action of a hyperbolic group, which is considered as a
version of Livschitz' theorem for hyperbolic dynamical systems.

\bigskip
\section{Preliminary}

We first recall the notions of Martin and Poisson boundaries, see
e.g.~\cite{r} and~\cite{w} for details. Let $\Gamma$ be a discrete
group with a symmetric finite set $S$ of generators. We denote
by~$|g|$ the word length and by $d(x,y)=|x^{-1}y|$ the word metric
with respect to $S$. Let $\mu$ be a finitely supported probability
measure on $\Gamma$. We shall assume that $\mu$ is nondegenerate
in the sense that the semigroup generated by the support of $\mu$
coincides with $\Gamma$. The measure~$\mu$ defines a random walk
on $\Gamma$ with transition probabilities $p(x,y)=\mu(x^{-1}y)$.
Let $P$ be the Markov operator with kernel $\{p(x,y)\}_{x,y}$.
Assume that the random walk is transient, that is, the Green
function
$$
G(x,y)=\sum^\infty_{n=0}p^{(n)}(x,y)
$$
is finite, where $\{p^{(n)}(x,y)\}_{x,y}$ is the kernel of $P^n$.
This is the case when $\Gamma$ is nonamenable. Moreover, in the
latter case there exist $C>0$ and $\lambda>0$ such that
\begin{equation} \label{1eExp}
G(x,y)\le Ce^{-\lambda d(x,y)}.
\end{equation}
Indeed, without loss of generality we may assume that the support
of $\mu$ is contained in $S$. By~\cite{dg} the spectral radius of
$P$ considered as an operator on $\ell^2(\Gamma)$ is less than
$1$. Hence there exists $z>1$ such that the operator
$$
G_z=\sum^\infty_{n=0}z^nP^n
$$
is bounded. In particular, there exists $c>0$ such that
$(G_z\delta_y,\delta_x)\le c$. Since by assumption we have a
nearest neighbour random walk, $p^{(n)}(x,y)=0$ for $n<d(x,y)$.
Hence
$$
G(x,y)\le z^{-d(x,y)}\sum^\infty_{n=0}z^np^{(n)}(x,y)
=z^{-d(x,y)}(G_z\delta_y,\delta_x)\le z^{-d(x,y)}c,
$$
which implies~(\ref{1eExp}).

The Martin kernel is defined by
$$
K(x,y)=\frac{G(x,y)}{G(e,y)},
$$
where $e\in\Gamma$ is the unit element. The Martin
compactification $\bar\Gamma$ of $\Gamma$ is the smallest
compactification such that $\Gamma\subset\bar\Gamma$ is discrete
and the functions $K(x,\cdot)$, $x\in\Gamma$, extend to continuous
functions on $\bar\Gamma$. The Martin boundary is
$\partial_M\Gamma=\bar\Gamma\backslash\Gamma$. The left action of
$\Gamma$ on itself extends to a continuous action on $\bar\Gamma$.

Let $\Omega=\prod^\infty_{n=0}\Gamma$ be the path space of our
random walk. For any point $g\in\Gamma$ we have a Markov measure
$\0P_g$ defined on paths starting at $g$. For $\0P_g$-a.e. path
$\underline{x}=\{x_n\}_n\in\Omega$ the sequence $\{x_n\}_n$
converges to a point on the boundary, so we get a map
$\Omega\to\partial_M\Gamma$. Denote by $\nu_g$ the image
of~$\0P_g$ under this map. The measure space
$(\partial_M\Gamma,\nu_e)$ is called the Poisson boundary of
$\Gamma$. We shall write $\nu$ instead of $\nu_e$. The measures
$\nu_g$, $g\in\Gamma$, are equivalent, and we have
$$
K(g,\omega)=\frac{d\nu_g}{d\nu}(\omega)=\frac{dg\nu}{d\nu}(\omega),
$$
where $g\nu$ is the measure defined by $g\nu(X)=\nu(g^{-1}X)$.

\smallskip

We next recall basic facts about hyperbolic groups, see
e.g.~\cite{gh}. As above we assume that $\Gamma$ is a finitely
generated group with a symmetric finite set $S$ of generators.
Denote the closed ball with center $x$ and radius $r$ by
$$
B(x,r)=\{g\in\Gamma\mid d(x,g)\leq r\}.
$$
For a subset $\Delta\subset\Gamma$, we write
$$N(\Delta,r)=\{g\in\Gamma\mid d(g,\Delta)\leq r\}$$
and $d\Delta=\{g\in\Gamma\mid d(g,\Delta)=1\}.$

The Gromov product is defined by the formula
$$(x|y)_z=\frac{1}{2}(d(x,z)+d(y,z)-d(x,y))$$
for $x,y,z\in\Gamma$. When $z$ is the unit element $e$, we simply
write $(x|y)=(x|y)_e$.

Let $\delta\geq 0$. The group $\Gamma$ is said to be
$\delta$-hyperbolic if
$$(x|y)\geq\min\{(x|z),(y|z)\}-\delta$$
for every $x,y,z\in\Gamma$.

If $\Gamma$ is $\delta$-hyperbolic, then every geodesic triangle
$\triangle=\{\alpha,\beta,\gamma\}$ in $\Gamma$ is $4\delta$-slim,
i.e.,
$$
\alpha\subset N(\beta\cup\gamma,4\delta),\ \ \beta\subset
N(\gamma\cup\alpha,4\delta),\ \ \gamma\subset
N(\alpha\cup\beta,4\delta).
$$

Given a geodesic triangle $\triangle=\{\alpha,\beta,\gamma\}$, we
associate to it in a natural manner a tripod, which is denoted by
$T_\triangle$, and a map $f_\triangle:\triangle\to T_\triangle$
whose restriction to each geodesic segment of $\triangle$ is
isometric. Then it is also known that if $\Gamma$ is
$\delta$-hyperbolic, then every geodesic triangle $\triangle$ is
$4\delta$-thin, i.e., if $f_\triangle(x)=f_\triangle(y)$, then we
have $d(x,y)\leq 4\delta$.

A sequence $\{x_i\}_{i\geq 1}$ in $\Gamma$ is said to converge to
infinity if $\lim_{i,j\to\infty}(x_i|x_j)=\infty$. The Gromov
boundary $\partial\Gamma$ is defined as the set of equivalence
classes of sequences converging to infinity in $\Gamma$, where two
sequences $\{x_i\}_{i\geq 1}$ and $\{y_i\}_{i\geq 1}$ are said to be
equivalent if $\lim_{i,j\to\infty}(x_i|y_j)=\infty$.

The Gromov product $(p|q)$ for $p,q\in\Gamma\cup\partial\Gamma$ is
defined by
$$(p|q)=\sup\liminf_{i,j\to\infty}(x_i|y_j)$$
where the $\sup$ above runs over all sequences $\{x_i\}_{i\geq 1}$
converging to $p$ and $\{y_i\}_{i\geq 1}$ converging to~$q$.
By $\delta$-hyperbolicity,
$$(p|q)-2\delta\leq \liminf_{i,j\to\infty}(x_i|y_j)\leq (p|q)$$
holds for all such $\{x_i\}_{i\geq 1}$ and $\{y_i\}_{i\geq 1}$.
Recall that $\Gamma\cup\partial\Gamma$ is compact equipped with
the base $\{B(x,r)\}\cup\{V_r(p)\}$, where
$$V_r(p)=\{q\in\Gamma\cup\partial\Gamma\mid (q|p)>r\}.$$
One can introduce a metric $\rho$ on $\partial\Gamma$ such that
$$
a^{-(\xi|\zeta)-c}\le \rho(\xi,\zeta)\le a^{-(\xi|\zeta)+c}
$$
for some $a>1$ and $c\ge0$. We call any such metric visual.

For $p,q\in\Gamma\cup\partial\Gamma$, we denote by $[p,q]$ the set
of all geodesic segments (or rays, lines) between $p$ and $q$.
For $p,q\in \partial\Gamma$ and any geodesic rays $x\in [e,p]$ and
$y\in [e,q]$, the quantity $(x(m)|y(n))$ is increasing both in $m$ and $n$.

Recall next that every geodesic triangle with $k$ vertices in
$\partial\Gamma$ and $3-k$ vertices in $\Gamma$ is
$4(k+1)\delta$-slim.

If $\{x(n)\}_{n=0}^\infty$ and $\{y(n)\}_{n=0}^\infty$ are
geodesic rays converging to the same point in $\partial\Gamma$,
then for all $n\geq 0$ we have
$$d(x(n),y(n))\leq d(x(0),y(0))+8\delta.$$
In addition, if $x(0)=y(0)$, then $d(x(n),y(n))\leq 4\delta$ for
$n\geq 0$.

\smallskip

Throughout the paper we assume that $\Gamma$ is a nonelementary
hyperbolic group (that is, it does not have a cyclic subgroup of
finite index) and study a random walk on it defined by a finitely
supported nondegenerate probability measure $\mu$. Then the group
is nonamenable, and its Martin boundary coincides with the Gromov
boundary by a result of Ancona~\cite{a}. We fix a finite symmetric
set $S$ of generators containing the support of $\mu$.

\bigskip


\section{Multiplicativity of the Green function along geodesic
segments}

For any points $x$, $y$ and $z$ we have $F(x,z)G(z,y)\le G(x,y)$,
where
$$
F(x,z)=\frac{G(x,z)}{G(z,z)}
$$
is the probability that a path starting at $x$ hits $z$. The main
technical result of Ancona~\cite{a} needed to identify the Martin
boundary of a hyperbolic group with its Gromov boundary is that up
to a factor the converse inequality is also true if $z$ lies on a
geodesic segment $\alpha\in[x,y]$. We shall need a slightly
stronger result saying that the same is true for the restriction
of the random walk to any subset containing a sufficiently large
neighbourhood of the segment. This is essentially contained
in~\cite{a}, but we shall give a detailed proof for completeness.

For a subset $\Delta\subset\Gamma$ consider the induced random
walk on $\Delta$ (to be precise, to get a random walk we have to
add a cemetery point to $\Delta$). We denote the corresponding
quantities using the subscript $\Delta$, so we write $P_\Delta$,
$G_\Delta$, and so on.

\bprop \label{Mult} There exist $R_0>0$  and $C\geq 1$ such that
if $x,y\in\Gamma$ and $v\in\Gamma$ lies on a geodesic segment
$\alpha\in[x,y]$, then
$$
G_\Delta(x,y)\leq C F_\Delta(x,v)G_\Delta(v,y)
$$
for any $\Delta\subset\Gamma$ containing $N(\alpha,R_0)$.
\eprop

For $\Lambda\subset\Gamma$ denote by $f^{\Lambda,(n)}(x,y)$ the
probability that a path from $x$ hits $\Lambda$ for the first time at the
$n$-th step and at $y\in\Lambda$, and by
$F^\Lambda(x,y)=\sum_nf^{\Lambda,(n)}(x,y)$ the probability that a
path from $x$ hits $\Lambda$ for the first time at $y$.
Dually, denote by $l^{\Lambda,(n)}(x,y)$ the probability that a path from
$x\in\Lambda$ stays in $\Gamma\backslash\Lambda$ till the $n$-th
step, when it passes through $y$, and put
$L^\Lambda(x,y)=\sum_nl^{\Lambda,(n)}(x,y)$.

Since $\mu$ is nondegenerate, there exists $K\in\mathbb{N}$ such
that $d(x,y)=1$ implies $p^{(k)}(x,y)>0$ for some $k\leq K$. Then
if $f$ is a superharmonic function, that is, $f\ge0$ and $Pf\le
f$, we have $p^{(k)}(x,y)f(y)\le f(x)$. Thus if we set
$$
C_1=\sup\{p^{(k)}(e,g)^{-1}\mid g\in S,\ k\le K,\ p^{(k)}(e,g)>0\},
$$
we get
\begin{equation} \label{eHarnack}
f(y)\le C_1f(x).
\end{equation}
The same is true for $x,y\in\Delta$, $d(x,y)=1$, if $f$ is
$P_\Delta$-superharmonic and $\Delta$ contains $N([x,y],K)$.
Moreover, similar inequalities with the same constant $C_1$ hold
for the random walk defined by the measure $\check\mu$,
$\check\mu(g)=\mu(g^{-1})$. For the random walk defined by
$\check\mu$ we write $\check P$, $\check P_\Delta$, and so on.

For $|z|<1/\rho(P)$, where $\rho(P)=\lim_{n\to\infty}p^{(n)}(e,e)^{1/n}<1$
is the spectral radius of our random walk, define
$$
G(x,y|z)=\sum^\infty_{n=0}p^{(n)}(x,y)z^n.
$$
We similarly introduce $F^\Lambda(x,y|z)$ and $L^\Lambda(x,y|z)$.

Observe that since $p^{(n)}_\Delta(x,y)\le p^{(n)}(x,y)$, the
functions $G_\Delta(x,y|z)$, $F^\Lambda_\Delta(x,y|z)$, and so on,
are well-defined for any $\Delta\subset\Gamma$ and
$\Lambda\subset\Delta$.

Fix a number $t$ such that $\rho(P)<t<1$, and set $s=1/t$.
Next choose $\ell\in\mathbb{N}$ such that $\ell\ge 8\delta$ and
$r\in\mathbb{N}$ such that $t^rC_1^{2\ell}\leq 1$. Finally put
$$
R_0=r+2\ell+K+1.
$$

For $x,y,v\in\Gamma$ define
$$
U_{x,v,y}=\{g\in\Gamma\mid (g|y)_x> d(x,v)\}.
$$
Since $(g|y)_x+(g|x)_y=d(x,y)$, it follows that if $v$ lies on a
geodesic segment between $x$ and $y$ then $\Gamma\setminus
U_{x,v,y}\subset U_{y,v,x}$.

\blem \label{MultL} There is $C'\geq 1$ such that
$$
G_\Delta(x,w)\leq C'F_\Delta(x,v)G_\Delta(v,w|s)
$$
for $w\in\Delta\cap N(U_{x,v,y},1)$, whenever $x,y\in\Gamma$ and
$v\in\Gamma$ lies on a geodesic segment $\alpha\in[x,y]$ with
$N(\alpha,R_0)\subset\Delta$. \elem

\bpf The proof follows that of \cite[Proposition 27.8]{w}. Let $m$
be the integer part of $d(x,v)/\ell$. Consider
$v_0,\dots,v_m\in\alpha$ between $x$ and $v$ such that $d(v_k,v)=(m-k)\ell$.
We denote $U_{x,v_k,y}$ by $W_k$. Then we have
\begin{enumerate}
\item[(i)] $v_k\in N(W_k,1)$ and $W_k\subset W_{k-1}$,
\item[(ii)] if $w\in W_k$
with $d(w,v_k)>r+1+8\delta$, then $B(w,r+1)\subset W_{k-1}$.
\end{enumerate}
Indeed, (i) is trivial. Let $w\in W_k$ with
$d(w,v_k)>r+1+8\delta$. Since $\triangle xwy$ is $4\delta$-thin
there is $v_k'$, which lies on a geodesic segment between $x$ and
$w$, such that $d(x,v_k')=d(x,v_k)$ and $d(v_k,v_k')\leq 4\delta$.
Let $w'\in B(w,r+1)$. If $(w'|w)_x\leq d(x,v'_k)$, then since
$\triangle xww'$ is $4\delta$-thin, there is $v_k''$, which lies
on a geodesic segment between $w'$ and $w$, such that
$d(v_k'',v_k')\leq 4\delta$. Then
$$
d(w,v_k)\le d(w,v_k'')+d(v_k'',v_k')+d(v_k',v_k)\le r+1+8\delta,
$$
which is a contradiction. Hence $(w'|w)_x> d(x,v'_k)=d(x,v_k)$.
Therefore
$$
(w'|y)_x\ge\min\{(w'|w)_x,(w|y)_x\} -\delta> d(x,v_k)-\delta\ge
d(x,v_{k-1}),
$$
so that $w'\in W_{k-1}$, and (ii) is proved.

The constant in the statement of the lemma will be
$C'=C_1^{2r+4\ell+2}$. By induction on $k$ we will show that
\begin{equation} \label{eMult1}
G_\Delta(x,w)\leq C'F_\Delta(x,v_k)G_\Delta(v_k,w|s)
\end{equation}
for $w\in\Delta\cap N(W_k,1)$.

For the case $k=0$, let $w\in\Delta$. Then $d(x,v_0)<\ell$. Since
$\Delta$ contains the $K$-neighbourhood of the geodesic segment
between $x$ and $v_0$, by the Harnack inequality~(\ref{eHarnack})
for $G_\Delta(\cdot,w)$ we obtain
$$
G_\Delta(x,w)\leq C_1^\ell G_\Delta(v_0,w)
\leq C_1^\ell G_\Delta(v_0,w|s).
$$
Similarly we also get
$$G_\Delta(v_0,v_0)\leq C_1^\ell G_\Delta(x,v_0)$$
and so
$$1\leq C_1^\ell F_\Delta(x,v_0).$$
It follows that (\ref{eMult1}) holds for $k=0$ and all
$w\in\Delta$.

Now suppose that inequality (\ref{eMult1}) holds for $k-1$. Using
again the Harnack inequality, we have
$$G_\Delta(v_{k-1},g|s)\leq C_1^\ell G_\Delta(v_k,g|s)$$
and
$$1\leq C_1^\ell F_\Delta(v_{k-1},v_k).$$
Therefore it follows from the induction hypothesis that for any
$g\in\Delta\cap W_{k-1}$
\begin{eqnarray*}
G_\Delta(x,g)&\leq& C'F_\Delta(x,v_{k-1})G_\Delta(v_{k-1},g|s) \\
&\leq&C'C_1^{2\ell}F_\Delta(x,v_{k-1})
F_\Delta(v_{k-1},v_k)G_\Delta(v_k,g|s) \\
&\leq&C'C_1^{2\ell}F_\Delta(x,v_k)G_\Delta(v_k,g|s).
\end{eqnarray*}

Let $w\in\Delta\cap N(W_k,1)$. Assume first that $d(w,v_k)\ge
r+2\ell+2$. Then $d(x,w)>r$. Indeed, let $w'\in W_k$ be such that
$d(w,w')\le1$. Then $d(w',v_k)\ge r+2\ell+1$, whence
$B(w',r+1)\subset W_{k-1}$ and also $B(w,r)\subset W_{k-1}$.
This proves the claim since $x\notin W_{k-1}$.

It follows that any path from $x$ to $w$ must pass through the set
$$
\Lambda=\{g\in\Delta\mid d(g,w)=r\}.
$$
Therefore using that $\Lambda\subset W_{k-1}$ and
$l^{\Lambda,(n)}_\Delta(g,w)=0$ for $n<r$ we get
\begin{eqnarray*}
G_\Delta(x,w)
&=&\sum_{g\in\Lambda}G_\Delta(x,g)L_\Delta^\Lambda(g,w) \\
&\leq&C'C_1^{2\ell}\sum_{g\in\Lambda}
F_\Delta(x,v_k)G_\Delta(v_k,g|s)t^rL_\Delta^\Lambda(g,w|s) \\
&\leq&C'F_\Delta(x,v_k)\sum_{g\in\Lambda}
G_\Delta(v_k,g|s)L_\Delta^\Lambda(g,w|s) \\
&=&C'F_\Delta(x,v_k)G_\Delta(v_k,w|s).
\end{eqnarray*}

Assume next that $d(v_k,w)\leq r+2\ell+1$. Then $\Delta$ contains
the $K$-neighbourhood of any geodesic segment between $v_k$ and
$w$. Hence by applying the Harnack inequality~(\ref{eHarnack}) to
the $\check P_\Delta$-superharmonic functions $G_\Delta(x,\cdot)$
and $G_\Delta(v_k,\cdot\,|s)$ we get
\begin{eqnarray*}
G_\Delta(x,w)&\leq&C_1^{r+2\ell+1}G_\Delta(x,v_k) \\
&\leq&C_1^{r+2\ell+1}F_\Delta(x,v_k)G_\Delta(v_k,v_k|s) \\
&\leq&C_1^{2r+4\ell+2}F_\Delta(x,v_k)G_\Delta(v_k,w|s).
\end{eqnarray*}
The proof is complete.
\epf

Applying the previous lemma to $\check\mu$ and using that $\check
F(g,h)=L(h,g)$, we get the following.

\bcor \label{MultC} We have
$$G_\Delta(w,y)\leq C'G_\Delta(w,v|s)L_\Delta(v,y)$$
for $w\in\Delta\cap N(U_{y,v,x},1)$. \ecor

We are now ready to prove the proposition.

\bpff{Proof of Proposition~\ref{Mult}} We follow the proof
of~\cite[Theorem~27.12]{w}. For $v=x,y$ the result is obvious.
Assume $v\ne x,y$. Then $x\not\in U_{x,v,y}$ and $y\in U_{x,v,y}$.
Since we have a nearest neighbour random walk, any path from $x$
to $y$ has to pass through the set $\Lambda=\Delta\cap dU_{x,v,y}$
on the way from $x$ to $y$. Therefore by Lemma~\ref{MultL} we get
\begin{equation} \label{eMult2}
G_\Delta(x,y)
=\sum_{w\in\Lambda}G_\Delta(x,w)L_\Delta^\Lambda(w,y)\leq
C'F_\Delta(x,v)
\sum_{w\in\Lambda}G_\Delta(v,w|s)L_\Delta^\Lambda(w,y).
\end{equation}
Moreover, for each $w\in\Lambda$ we have $w\in
N(U_{y,v,x},1)\cap\Delta$. Hence by Corollary~\ref{MultC} we also have
$$
G_\Delta(w,y)\leq C'G_\Delta(w,v|s)L_\Delta(v,y).
$$
Recall, see e.g.~\cite[Lemma~27.5]{w}, that if $f$ is a
superharmonic function then the minimal superharmonic function
dominating $f$ on $\Lambda$ is
$$
f^\Lambda(g)=\sum_{h\in\Lambda}F^\Lambda(g,h)f(h).
$$
Applying this to $G_\Delta(\cdot,v|s)$ we conclude that for all
$a\in \Delta$
$$
\sum_{w\in\Lambda}F^\Lambda_\Delta(a,w)G_\Delta(w,y)\le
C'\sum_{w\in \Lambda}F^\Lambda_\Delta(a,w) G_\Delta(w,v|s)L_\Delta(v,y)
\le C'G_\Delta(a,v|s)L_\Delta(v,y),
$$
so that
\begin{equation} \label{eMult3}
\sum_{w\in\Lambda}G_\Delta(a,w)L_\Delta^\Lambda(w,y)=
\sum_{w\in\Lambda}F^\Lambda_\Delta(a,w)G_\Delta(w,y) \leq
C'G_\Delta(a,v|s)L_\Delta(v,y).
\end{equation}

Combing (\ref{eMult2}) and (\ref{eMult3}) with the resolvent
equation
$$
sG_\Delta(v,w|s)-G_\Delta(v,w)
=(s-1)\sum_{a\in \Delta}G_\Delta(v,a|s)G_\Delta(a,w),
$$
we get
\begin{eqnarray*}
G_\Delta(x,y)&\leq&C'F_\Delta(x,v)
\sum_{w\in\Lambda}G_\Delta(v,w|s)L_\Delta^\Lambda(w,y) \\
&\leq&C'F_\Delta(x,v)\sum_{w\in\Lambda}
\left\{\frac{1}{s}G_\Delta(v,w)+\left(1-\frac{1}{s}\right)
\sum_{a\in \Delta}G_\Delta(v,a|s)G_\Delta(a,w)\right\}
L_\Delta^\Lambda(w,y) \\
&\leq&{C'}^2F_\Delta(x,v)\left\{\frac{1}{s}G_\Delta(v,v|s)
+\left(1-\frac{1}{s}\right)\sum_{a\in
\Delta}G_\Delta(v,a|s)G_\Delta(a,v|s)\right\}L_\Delta(v,y).
\end{eqnarray*}
Choose $s'$ such that $s<s'<1/\rho(P)$. Then using the resolvent
equation again and the fact
$$p_\Delta^{(n)}(v,v)\leq p^{(n)}(v,v)\leq\rho(P)^n,$$
we have
\begin{align*}
\frac{1}{s}G_\Delta(v,v|s)&+\left(1-\frac{1}{s}\right)
\sum_{a\in\Delta}G_\Delta(v,a|s)G_\Delta(a,v|s) \\
&\leq\frac{1}{s}G_\Delta(v,v|s)+\left(1-\frac{1}{s}\right)
\sum_{a\in\Delta}G_\Delta(v,a|s')G_\Delta(a,v|s) \\
&=\frac{1}{s}G_\Delta(v,v|s)+\left(1-\frac{1}{s}\right)
\frac{1}{s'-s}\{s'G_\Delta(v,v|s')-sG_\Delta(v,v|s)\} \\
&\leq\frac{1}{s}G_\Delta(v,v|s)+\left(1-\frac{1}{s}\right)
\frac{s'}{s'-s}G_\Delta(v,v|s') \\
&\leq\frac{1}{s(1-\rho(P)s)}+\left(1-\frac{1}{s}\right)
\frac{s'}{(s'-s)(1-\rho(P)s')}.
\end{align*}
Denote the last expression by $C''$. Then, since
$G_\Delta(v,v)\geq 1$, if we set $C={C'}^2C''$, we obtain
$$
G_\Delta(x,y)\leq CF_\Delta(x,v)G_\Delta(v,v)L_\Delta(v,y)
=CF_\Delta(x,v)G_\Delta(v,y),
$$
and the proposition is proved. \epf

It will be convenient to have a version of the above proposition
for the case when $v$ is only close to a geodesic segment.

\bcor \label{MultC2} For any $R_1$ there exists $C'\ge1$ such that
if $x,y\in\Gamma$, $\alpha\in[x,y]$ and $v\in N(\alpha,R_1)$, then
$$
G_\Delta(x,y)\leq C' F_\Delta(x,v)G_\Delta(v,y)
$$
for any $\Delta\subset\Gamma$ containing $N(\alpha,R_0+R_1)$.
\ecor

\bpf By Proposition~\ref{Mult} applied to $v'\in\alpha$ such that
$d(v,v')\le R_1$ we have
$$
G_\Delta(x,y)\leq C G_\Delta(x,v')G_\Delta(v',y).
$$
Since $R_0\ge K$, the set $\Delta$ contains $N([v,v'],K)$. So
applying the Harnack inequality~(\ref{eHarnack}) twice we get
$$
G_\Delta(x,y)\le CC^{2R_1}_1G_\Delta(x,v)G_\Delta(v,y)=
CC^{2R_1}_1F_\Delta(x,v)G_\Delta(v,v)G_\Delta(v,y),
$$
so that we can take $C'=CC^{2R_1}_1G(e,e)$. \epf

\bigskip

\section{A Harnack inequality at infinity and the H\"older
condition}

In~\cite{l} Ledrappier proves that in the case of a free group the
Martin kernel is H\"older continuous, which is a discrete analogue
of a result of Anderson and Schoen~\cite{as}. Our goal in this
section is to extend this result to hyperbolic groups.

The first step is to prove an analogue of~\cite[Theorem 3.1]{l}, a
Harnack inequality at infinity.

We say that a function $u$ on $\Gamma$ is harmonic on
$\Lambda\subset\Gamma$ if $(Pu)(g)=u(g)$ for $g\in\Lambda$. Since
we consider a nearest neighbour random walk, any function
which coincides with $u$ on $N(\Lambda,1)$ is harmonic on
$\Lambda$.

For $r\geq 0$ and $p\in\Gamma$ we define
$$C_r(p)=\{g\in\Gamma\mid (g|e)_p>r\}.$$
We fix an integer $R$ such that $R\ge R_0+14\delta+2$, where $R_0$
is from Proposition~\ref{Mult}.

\bprop \label{Harnack1} There exists a constant $B\geq 1$
satisfying the following property. Let $m\in\mathbb{N}$,
$k\in\mathbb{N}$ with $0\leq k\leq m-1$, $p\in\Gamma$ with
$|p|=3Rm$ and $\alpha\in[p,e]$. Assume that $u,v$ are functions on
$\Gamma$ which are positive and harmonic on $C_{3Rk}(p)$ and
vanish at infinity on $C_{3Rk}(p)$. Then we have
$$\frac{u(g)}{u(\alpha(3Rk+2R))}\leq
B \frac{v(g)}{v(\alpha(3Rk+2R))}
\ \ \hbox{for}\ g\in C_{3R(k+1)}(p).$$ \eprop

\bpf Let $T(\underline{x})$ be the first time a path
$\underline{x}$ from $g\in C_{3Rk+1}(p)$ hits
$\Lambda=dC_{3Rk+1}(p)$, and put $T_n=\min\{T,n\}$. Note that
$\Lambda\subset C_{3Rk}(p)$. Since any path from $g$ stays in
$C_{3Rk+1}(p)$ till it hits~$\Lambda$ (if ever), and $u$ is
harmonic on $C_{3Rk+1}(p)$, we have
$$
u(g)=\int_\Omega u(x_{T_n})d\0P_g(\underline{x}).
$$
Since our random walk is transient, almost every path either hits
$\Lambda$, or goes to infinity. Since $u$ vanishes at infinity on
$C_{3Rk}(p)$, letting $n\to\infty$ we thus get
$$
u(g)=\int_{\{T<\infty\}}u(x_{T})d\0P_g(\underline{x})
=\sum_{h\in\Lambda}F^\Lambda(g,h)u(h).
$$

If $g\in C_{3Rk+R}(p)$, then any path from $g$ to $h\in\Lambda$
passes through $\Theta=dC_{3Rk+R}(p)$. Hence, denoting
$C_{3Rk+1}(p)$ by $\Delta$, we can write
\begin{equation} \label{eHarInf1}
F^\Lambda(g,h)=\sum_{a\in\Theta}\sum_{b\in\Delta\backslash\Theta}
G_\Delta(g,a)L^\Theta_\Delta(a,b)p(b,h).
\end{equation}

If $g\in C_{3R(k+1)}(p)$ and $a\in dC_{3Rk+R}(p)$ then
$(a|e)_p>3Rk+R-1$, whence by $\delta$-hyperbolicity
$$
(a|g)_p\ge\min\{(a|e)_p,(e|g)_p\}-\delta>3Rk+R-1-\delta.
$$
Since $\triangle pga$ is $4\delta$-thin, we conclude that for any
$b$ lying on a geodesic segment $\gamma\in[g,a]$ we have
$$
(b|g)_p>3Rk+R-1-3\delta.
$$
Therefore
$$
(b|e)_p\ge\min\{(b|g)_p,(g|e)_p\}-\delta>3Rk+R-1-4\delta.
$$
It follows that $N(\gamma,R_0+8\delta)$ is contained in
$C_{3Rk+1}(p)$. On the other hand,
$$
3Rk+R\geq(a|e)_p\geq\min\{(a|g)_p,(g|e)_p\}-\delta=(a|g)_p-\delta,
$$
and so $(a|g)_p<3Rk+2R$. Hence if $\beta\in[p,g]$ then
$\beta(3Rk+2R)\in N(\gamma,4\delta)$. Denote $\alpha(3Rk+2R)$ by
$g_0$. Then we also have $d(\beta(3Rk+2R),g_0)\leq 4\delta$,
because $(g|e)_p>3R(k+1)$. It follows that $g_0\in
N(\gamma,8\delta)$. Hence by Corollary~\ref{MultC2} there exists
$B\ge1$ such that
$$
F_\Delta(g,g_0) G_\Delta(g_0,a)\le G_\Delta(g,a)\le B
F_\Delta(g,g_0) G_\Delta(g_0,a)
$$
for any $g\in C_{3R(k+1)}(p)$ and $a\in\Theta$. By virtue
of~(\ref{eHarInf1}) we get
$$
F_\Delta(g,g_0)F^\Lambda(g_0,h)\le F^\Lambda(g,h)\le B
F_\Delta(g,g_0)F^\Lambda(g_0,h).
$$
It follows that
$$
F_\Delta(g,g_0)u(g_0)\le u(g)\le BF_\Delta(g,g_0)u(g_0).
$$
Since the same inequalities hold for $v$, we get the result. \epf

Next we shall prove an analogue of \cite[Lemma 3.10]{l}.

\blem \label{Harnack2} There exist $B'\ge1$ and $0\le\sigma<1$
such that for any $m\in\mathbb{N}$, any pair
$\xi,\eta\in\partial\Gamma$ with $(\xi|\eta)>3R(m+2)$ and any
geodesic ray $\{x(n)\}_{n=0}^\infty$ $\in[e,\xi]$
the function
$$\varphi(g)=\frac{K(g,\xi)}{K(g,\eta)}$$
has the property
$$|\varphi(g)-\varphi(h)|\leq B'\sigma^{k-1}$$
for $g,h\in C_{3R(k+1)}(p)$ with $1\leq k\leq m$, where
$p=x(3R(m+2))$. \elem

\bpf The functions $K(\cdot,\xi)$ and $K(\cdot,\eta)$ are positive
and harmonic. We claim that they vanish at infinity on
$C_{3R}(p)$. Indeed, choose a geodesic ray
$\{y(n)\}_n\in[e,\eta]$. Since
$$
\lim_n(x(n)|y(n))\ge(\xi|\eta)-2\delta
$$
by $\delta$-hyperbolicity, we have $(x(n)|y(n))>3R(m+2)-2\delta$
for any sufficiently large $n$. Since the triangle $\triangle
ey(n)x(n)$ is $4\delta$-thin, we then get
\begin{equation} \label{eHarnack2}
(y(n)|p)>3R(m+2)-4\delta.
\end{equation}
Let $g\in C_{3R}(p)$. Choose $\alpha\in[p,g]$, $\beta\in[p,y(n)]$,
$\gamma\in[g,y(n)]$. Denote $x(3R(m+1))$ by $q$. Since
$(g|e)_p>3R$, using $4\delta$-thinness of $\triangle gep$ we can
find $a\in\alpha$ such that $d(a,q)\le4\delta$. Since
$\{\alpha,\beta,\gamma\}$ is $4\delta$-slim, we can find
$b\in\beta\cup\gamma$ such that $d(a,b)\le4\delta$. If $b\in\beta$
then using (\ref{eHarnack2}) and that $\triangle ey(n)p$ is
$4\delta$-thin we get
$$
(b|p)> 3R(m+2)-6\delta,
$$
whence
$$
(q|p)> 3R(m+2)-14\delta.
$$
This is a contradiction as $(q|p)=3R(m+1)$. Hence $b\in\gamma$. It
follows that $q\in N(\gamma,8\delta)$. By Corollary~\ref{MultC2}
there exists a constant $C$ such that
$$
G(g,y(n))\le CG(g,q)G(q,y(n)),
$$
so that
$$
K(g,\eta)=\lim_{n\to\infty}\frac{G(g,y(n))}{G(e,y(n))}\le
\lim_{n\to\infty}
C\frac{G(g,q)G(q,y(n))}{F(e,q)G(q,y(n))}=C\frac{G(g,q)}{F(e,q)}.
$$
Since $G(\cdot,q)$ vanishes at infinity by (\ref{1eExp}), we
conclude that $K(\cdot,\eta)$ vanishes at infinity on $C_{3R}(p)$.
The same is true for $K(\cdot,\xi)$ as we could take $\eta=\xi$.

\smallskip

By Proposition~\ref{Harnack1} we have
\begin{equation} \label{eHarnack3}
\frac{K(g,\xi)}{K(x(3Rm+R),\xi)}
\leq B\frac{K(g,\eta)}{K(x(3Rm+R),\eta)}
\ \ \mbox{for}\ g\in C_{6R}(p).
\end{equation}
Let $\{y(n)\}^\infty_{n=0}\in[e,\eta]$. Since as above
$(x(n)|y(n))>3R(m+2)-2\delta$ for any $n$ sufficiently large, we
have
$$
d(x(3Rm+R),y(3Rm+R))\le4\delta.
$$
Denote $x(3Rm+R)$ by $a$ and $y(3Rm+R)$ by $b$. We have
$$
K(a,\xi)=\lim_n\frac{G(a,x(n))}{G(e,x(n))}\le
\lim_n\frac{G(a,x(n))}{F(e,a)G(a,x(n))}=\frac{1}{F(e,a)}
=\frac{G(e,e)}{G(e,a)}.
$$
On the other hand, by Corollary~\ref{MultC2} there exists a
constant $C$ such that
$$
G(e,y(n))\le CG(e,b)G(b,y(n)),
$$
whence
$$
K(b,\eta)\ge\lim_n\frac{G(b,y(n))}{CG(e,b)G(b,y(n))}
=\frac{1}{CG(e,b)}.
$$
Since $G(e,b)\le C_1^{4\delta}G(e,a)$ by (\ref{eHarnack}), we thus
get
$$
\frac{K(a,\xi)}{K(b,\eta)}\le CC_1^{4\delta}G(e,e).
$$
Put $B'=B CC_1^{4\delta}G(e,e)$. Then in view of (\ref{eHarnack3})
we have
$$
0\leq\varphi(g)=\frac{K(g,\xi)}{K(g,\eta)}\leq B' \ \
\mbox{for}\ g\in C_{6R}(p).
$$

Denote $K(g,\xi)$ by $u(g)$ and $K(g,\eta)$ by $v(g)$. For $2\leq
k\leq m+1$, define
$$
\underline{\varphi}_k=\inf_{g\in
C_{3Rk}(p)}\varphi(g),\ \ \ \ \ \overline{\varphi}_k=\sup_{g\in
C_{3Rk}(p)}\varphi(g)
$$
and $u_k(g)=u(g)-\underline{\varphi}_kv(g)$. Note that $u_k$ is
positive on $C_{3Rk}(p)$. Applying Proposition~\ref{Harnack1} we
have
$$
\frac{D_k}{B}\leq\frac{u_k(g)}{v(g)}\leq D_kB\ \ \mbox{for}\ g\in
C_{3R(k+1)}(p),
$$
where
$$
D_k=\frac{u_k(x(3R(k+1)+R))}{v(x(3R(k+1)+R))}.
$$
We therefore obtain
$$\sup_{g\in C_{3R(k+1)}}\frac{u_k(g)}{v(g)}
\leq B^2\inf_{g\in C_{3R(k+1)}}\frac{u_k(g)}{v(g)},
$$
and so
$$\overline{\varphi}_{k+1}-\underline{\varphi}_k
\leq B^2(\underline{\varphi}_{k+1}-\underline{\varphi}_k).
$$
Similarly by using $v_k(g)=\overline{\varphi}_kv(g)-u(g)$, we also
get
$$\overline{\varphi}_k-\underline{\varphi}_{k+1}\leq
B^2(\overline{\varphi}_k-\overline{\varphi}_{k+1}).$$
Hence if we set
$\omega_k=\overline{\varphi}_k-\underline{\varphi}_k$, then we
have
$$\omega_{k+1}+\omega_k\leq B^2(\omega_k-\omega_{k+1}),$$
and thus $\omega_{k+1}\leq\sigma\omega_k$, where
$$0\le\sigma=\frac{B^2-1}{B^2+1}<1.$$
Therefore for $g,h\in C_{3R(k+1)}(p)$, we have
$$
|\varphi(g)-\varphi(h)|\leq\omega_{k+1}\leq\cdots
\leq\sigma^{k-1}\omega_2\leq B'\sigma^{k-1},
$$
and the proof of the lemma is complete. \epf

We can now prove that the Martin kernel is H\"older continuous.

\bthm \label{Holder} There exist $0\le\tau<1$ and for any
$g\in\Gamma$ a constant $H_g\ge0$ such that for
$\xi,\eta\in\partial\Gamma$,
$$|K(g,\xi)-K(g,\eta)|\leq H_g\tau^{(\xi|\eta)}.$$
\ethm

\bpf
The proof repeats that of~\cite[Theorem 3.7]{l}.
Let $g\in\Gamma$. Let $n\in\0N$ be such that
$$
9R+|g|\ge3Rn>6R+|g|,
$$
where $R$ is from Proposition~\ref{Harnack1}.
Let $\xi,\eta\in\partial\Gamma$ be such that $(\xi|\eta)>3Rn$.
Then
$$
3R(m+3)\geq(\xi|\eta)> 3R(m+2)
$$
for some $m\ge n-2$. Let $\{x(n)\}^\infty_{n=0}\in[e,\xi]$. Put
$p=x(3R(m+2))$. Then
$$
(g|e)_p\ge|p|-|g|\ge 3Rn-|g|>6R,
$$
so that $g\in C_{6R}(p)$. Take $k\in\mathbb{N}$ such that
$3R(k+2)\geq(g|e)_p>3R(k+1)$. By applying Lemma~\ref{Harnack2} to
$g, e\in C_{3R(k+1)}(p)$, we have
$$\left|\frac{K(g,\xi)}{K(g,\eta)}-1\right|\leq B'\sigma^{k-1}.$$
Hence
$$|K(g,\xi)-K(g,\eta)|\leq K(g,\eta)B'\sigma^{k-1},$$
Note that
\begin{eqnarray*}
(\xi|\eta)&\leq& 3R(m+3)
=|p|+3R=(g|e)_p+(g|p)_e+3R \\
&\leq&3R(k+2)+|g|+3R=3R(k-1)+|g|+12R.
\end{eqnarray*}
Thus if we put $\tau=\sigma^{1/3R}$
and $M_g=\sup\{K(g,\zeta)\mid \zeta\in\partial\Gamma\}$, then
\begin{equation} \label{eHolder}
|K(g,\xi)-K(g,\eta)|\leq B'K(g,\eta)\tau^{-|g|-12R}\tau^{(\xi|\eta)}
\le B'M_g\tau^{-|g|-12R}\tau^{(\xi|\eta)}.
\end{equation}
for any $\xi$ and $\eta$ such that $(\xi|\eta)>3Rn$. On the other
hand, if $(\xi|\eta)\le3Rn<12R+|g|$ then
$$
|K(g,\xi)-K(g,\eta)|\leq 2M_g
\le 2M_g\tau^{-|g|-12R}\tau^{(\xi|\eta)}.
$$
Thus it suffices to put $H_g=2B'M_g\tau^{-|g|-12R}$. \epf

From the proof of the above theorem we also get the following.

\bcor \label{HolderC} There exist $H\ge0$ and $L_0>0$ such that
$$\left|\frac{K(g,\xi)}{K(g,\eta)}-1\right|
\leq H\tau^{(\xi|\eta)-|g|}$$
whenever $(\xi|\eta)-|g|\geq L_0$. \ecor

\bpf Using the notation from the proof of the theorem, take
$L_0>9R\ge 3Rn-|g|$. Then if $(\xi|\eta)\ge L_0+|g|$, we have
$(\xi|\eta)
>3Rn$, so that by the first inequality in (\ref{eHolder})
$$
\left|\frac{K(g,\xi)}{K(g,\eta)}-1\right|
\leq B'\tau^{-12R}\tau^{(\xi|\eta)-|g|}.
$$
Thus we can take $H=B'\tau^{-12R}$. \epf

\bigskip

\section{A Gibbs-like property of a harmonic measure}

For $\xi\in\partial\Gamma$ and $R>0$, we define $U(\xi,R)$ to be
the set of all $\eta\in\partial\Gamma$ such that for any pair of
geodesic rays $\{x(n)\}_{n=0}^\infty\in[e,\xi]$ and
$\{y(n)\}_{n=0}^\infty\in[e,\eta]$, we have
$$\lim_{n\to\infty}(x(n)|y(n))>R.$$
Remark that the sequence $\{(x(n)|y(n))\}_{n=0}^\infty$ is
nondecreasing and thus the above limit always exists. These sets
are considered as hyperbolic versions of cylindric sets. Note that
if $\zeta\in U(\xi,R)$ then $(\zeta|\xi)>R$, and by
$\delta$-hyper\-bolicity if $(\zeta|\xi)>R$ then $\zeta\in
U(\xi,R-2\delta)$. We also have that
\begin{equation} \label{eGibbs}
\hbox{if}\ \ \eta\in U(\xi,R+\delta),\ \ \hbox{then}\ \
U(\xi,R+\delta)\subset U(\eta,R).
\end{equation}

\smallskip

The following property of the harmonic measure $\nu=\nu_e$ on
$\partial\Gamma$ reminds of a Gibbs measure.

\bthm\label{Gibbs} There exists $D\geq 1$ such that for every
$\xi\in\partial\Gamma$ and $\{x(n)\}_{n=0}^\infty\in[e,\xi]$, we
have
$$\frac{1}{D}\leq\frac{\nu(U(\xi,R))}{F(e,x(R))}\leq D \ \
\mbox{for}\ R\in\mathbb{N}.$$
\ethm

We need the following lemma to prove the theorem.

\blem \label{Gibbs1} Let $N\in\mathbb{N}$ be such that
$N>13\delta$. Then for any $m\in\mathbb{N}$,
$\xi\in\partial\Gamma$ and $\{x(n)\}_{n=0}^\infty\in[e,\xi]$, we
have
$$U(x(m)^{-1}\xi,N)\subset x(m)^{-1}U(\xi,m).$$
\elem

\bpf Let $\zeta\in x(m)U(x(m)^{-1}\xi,N)$ and
$\{y(n)\}^\infty_{n=0}\in[x(m),\zeta]$. Then
$\{x(m)^{-1}y(n)\}^\infty_{n=0}\in[e,x(m)^{-1}\zeta]$ and
$\{x(m)^{-1}x(m+n)\}^\infty_{n=0}\in[e,x(m)^{-1}\xi]$. Since
$x(m)^{-1}\zeta\in U(x(m)^{-1}\xi,N)$, we have
$$
(x(m)^{-1}y(n)|x(m)^{-1}x(m+n))>N
$$
for any $n$ sufficiently large, whence
$$
4\delta\ge d(x(m)^{-1}y(N),x(m)^{-1}x(m+N))=d(y(N),x(m+N)).
$$

Now let $\{z(n)\}^\infty_{n=0}\in[e,\zeta]$. Since the triangle
$x([0,m])\cup y\cup z$ is $8\delta$-slim, there exists $a\in
x([0,m])\cup z$ such that $d(y(N),a)\le8\delta$. Then
$d(a,x(m+N))\le12\delta$. Since $N>12\delta$, we cannot have $a\in
x([0,m])$. Hence $a\in z$. Let $a=z(k)$. As
$d(z(k),x(m+N))\le12\delta$, we have $|m+N-k|\le12\delta$, and
therefore
$$
d(z(m+N),x(m+N))\le24\delta.
$$
It follows that for any $n\ge N+m$ we have
$$
(z(n)|x(n))\ge(z(m+N)|x(m+N))\ge m+N-12\delta.
$$
Then, for any $\{x'(n)\}^\infty_{n=0}\in[e,\xi]$,
$$
\lim_n(z(n)|x'(n))\ge m+N-13\delta>m.
$$
Thus $\zeta\in U(\xi,m)$. \epf

\bpff{Proof of Theorem~\ref{Gibbs}} Let $\eta\in U(\xi,R)$ and
$\{y(n)\}_{n=0}^\infty\in[e,\eta]$. Then $(x(m)|y(m))>R$ for
sufficiently large $m$, whence $d(x(R),y(R))\le4\delta$. Thanks to
Corollary~\ref{MultC2} there exists a positive constant $C$ such
that
$$F(e,x(R))G(x(R),y(m))\leq G(e,y(m))\leq CF(e,x(R))G(x(R),y(m)).$$
Dividing by $G(e,y(m))$ and letting $m\to\infty$ we get
$$F(e,x(R))\frac{d\nu_{x(R)}}{d\nu}(\eta)\leq 1
\leq CF(e,x(R))\frac{d\nu_{x(R)}}{d\nu}(\eta).$$ Integrating over
$U(\xi,R)$ we obtain
$$\nu_{x(R)}(U(\xi,R))\leq\frac{\nu(U(\xi,R))}{F(e,x(R))}
\leq C\nu_{x(R)}(U(\xi,R))\leq C.$$ Fix $N>13\delta$. Since
$$
\nu_{x(R)}(U(\xi,R))=\nu(x(R)^{-1}U(\xi,R))\ge\nu(U(x(R)^{-1}\xi,N))
$$
by Lemma~\ref{Gibbs1}, to prove the theorem it suffices to show
that there exists a positive constant~$D'$ such that
$\nu(U(x(R)^{-1}\xi,N))\geq D'$ for any $\xi\in\partial\Gamma$,
$R\in\mathbb{N}$ and $\{x(n)\}_{n=0}^\infty\in[e,\xi]$.

Assume such a constant does not exist. Then there exist
$\xi_k\in\partial\Gamma$, $m_k\in\0N$ and
$\{x_k(n)\}_{n=0}^\infty\in[e,\xi_k]$ such that
$$
\nu(U(x_k(m_k)^{-1}\xi_k,N))\to0\ \ \hbox{as}\ \ k\to\infty.
$$
Since $\partial\Gamma$ is compact, by taking a subsequence if
necessary, we may assume that the sequence
$\{\eta_k=x_k(m_k)^{-1}\xi_k\}_k$ converges to a point
$\zeta\in\partial\Gamma$. Then~(\ref{eGibbs}) implies that for
sufficiently large~$k$ the set $U(\eta_k,N)$ contains $U(\zeta,
N+\delta)$. Hence $\nu(U(\zeta,N+\delta))=0$. This is a
contradiction, because the action of $\Gamma$ on $\partial\Gamma$
is minimal and hence any open set has positive measure. \epf

\bigskip

\section{A Livschitz type theorem}
Livschitz' theorem ~\cite[Theorem~19.2.1]{kh} says that every
H\"older continuous cocycle of a topologically transitive
hyperbolic dynamical system is a coboundary given by a H\"older
continuous function with the same exponent. In this section we
establish a Livschitz type theorem for the boundary action of a
hyperbolic group.

It is known that every infinite order element $g\in\Gamma$ acts on
$\partial\Gamma$ as a hyperbolic homeomorphism, i.e., there are
exactly two fixed points $g^+$ and $g^-$ in $\partial\Gamma$ such
that $g^+$ is stable and $g^-$ is unstable. For any open subsets
$U^\pm\subset\partial\Gamma$ with $g^+\in U^+$ and $g^-\in U^-$,
it holds that $g^n(\partial\Gamma\setminus U^-)\subset U^+$ for
sufficiently large $n\geq 0$. For $h\in\Gamma$, the sequence
$\{g^nh\}_{n=1}^\infty$ converges to $g^+$ and
$\{g^{-n}h\}_{n=1}^\infty$ converges to $g^{-}$.

\bthm\label{Liv} Let $G$ be a group with a two-sided invariant
metric $\rho$ such that $G$ is complete with respect to $\rho$. We
assume that $c\colon\Gamma\times \pG\rightarrow G$ is a H\"older
continuous cocycle, that is, \begin{itemize} \item[(1)] there
exist positive constants $0<\tau<1$ and $A_g$ for each $g\in
\Gamma$ satisfying
$$\rho(c(g,\xi),c(g,\eta))\leq A_g\tau^{(\xi|\eta)},\quad
\forall \xi,\eta\in \pG,$$ \item[(2)] the cocycle identity
$c(gh,\omega)=c(g,\omega)c(h,g^{-1} \omega)$ holds for every
$g,h\in \Gamma$ and $\omega\in\pG$.
\end{itemize}
We assume $c(g,g^{+})=c(g,g^{-})=e$ for every infinite order
element $g\in \Gamma$. Then there exists a continuous map
$b\colon\pG\rightarrow G$ satisfying
$$c(g,\xi)=b(\xi)b(g^{-1}\xi)^{-1},\quad \forall g\in \Gamma,\;
\forall \xi\in \pG.$$
The map $b$ satisfies H\"older continuity with the same exponent as
$c$, that is, there exists a positive constant $A$ such that
$$\rho(b(\xi),b(\eta))\leq A\tau^{(\xi|\eta)},\quad
\forall \xi,\eta\in \pG.$$
\ethm

Before proving the theorem, we show a few lemmas.

\blem \label{Livlem1}
Let $g\in \Gamma$ be an infinite order element and
$r>0$.
Then there exists a constant $C(g,r)$ depending only on
$g$, $r$ and $\delta$ such that the following holds for any
$\xi,\eta \in \pG\setminus V_r(g^{-})$ and any natural number $n$:
$$(g^n\xi|g^n\eta)\geq |g^{-n}|+(\xi|\eta)-C(g,r).$$
\elem

\bpf Let $x\in [e,\xi]$ and $y\in [e,\eta]$. Then
$$(g^nx(m)|g^ny(m))=(x(m)|y(m))+|g^{-n}|-(g^{-n}|x(m))-(g^{-n}|y(m))$$
and so
\begin{eqnarray*}
(g^n\xi|g^n\eta)&\geq& \lim_{m\to\infty}(g^nx(m)|g^ny(m))\\
&\geq& (\xi|\eta)-2\delta+|g^{-n}|
-\liminf_{m\to\infty}(g^{-n}|x(m))-\liminf_{m\to\infty}(g^{-n}|y(m)).
\end{eqnarray*}
Let $z\in [e,g^-]$. Since $\{z(m)\}_{m=0}^\infty$ and
$\{g^{-n}\}_{n=0}^\infty$ converge to the same point $g^-$, there
exists $m_0$ such that $(g^{-n}|z(m))>r+\delta$ holds for any
$m,n\geq m_0$. Since $\xi\notin V_r(\xi)$ and $\eta\notin
V_r(\eta)$, we have $(x(m)|z(m))\leq r$ and $(y(m)|z(m))\leq r$
for any $m\in\0 N$. Thus $\delta$-hyperbolicity implies
$$\min\{(g^{-n}|x(m)),(g^{-n}|z(m))\}-\delta\leq (x(m)|z(m))\leq r$$
and so $(g^{-n}|x(m))\leq r+\delta$ holds for any $m,n\geq m_0$.
In particular, for any $n\geq m_0$ we get
$$\liminf_{m\to\infty}(g^{-n}|x(m))\leq r+\delta.$$
Therefore
$$(g^n\xi|g^n\eta)\geq (\xi|\eta)+|g^{-n}|-2\delta
-2\max\{r+\delta,|g^{-k}|\colon 0\leq k<m_0\}$$ holds for any
natural number $n$. \epf

\blem \label{Livlem2} Let $g\in \Gamma$ be an infinite order
element and $h\in \Gamma$. We assume that $h g^{+}\neq g^{-}$.
Then there exists $m_{g,h}\in \0 N$ such that \begin{itemize}
\item [(1)] the element $g^nh$ is of infinite order for any $n\geq
m_{g,h}$, \item [(2)] the sequence
$\{(g^nh)^{+}\}_{n=m_{g,h}}^\infty$ converges to $g^+$.
\end{itemize} \elem

\bpf (1). Let $U$ be a neighbourhood of $g^+$ such that
$g^{-}\notin h\overline{U}$, where $\overline{U}$ is the closure
of $U$. Then there exists $m\in \0 N$ such that $g^nhU$ is
strictly included in $U$ for any $n\geq m$. For such~$n$, the
$g^nh$-orbit of any point in $U\setminus g^nhU$ is an infinite
set, which shows that $g^nh$ is of infinite order.

(2). The above argument shows that $\omega\ne(g^nh)^-$ for any
$\omega\in U\setminus g^nhU$. Hence the sequence
$\{(g^nh)^k\omega\}_{k=1}^\infty$ converges to $(g^nh)^+$, which
shows that $(g^nh)^+$ belongs to $\overline U$. Since this holds
for every $U$ as above and sufficiently large $n$, we get the
statement. \epf

\blem \label{Livlem3} Let $g\in \Gamma$ be an infinite order
element and $n\in \0 N$. For  $\xi\in \pG\setminus \{g^-\}$ we set
$b_{g,n}(\xi)=c(g^{-n},\xi)$. Then the sequence
$\{b_{g,n}\}_{n=1}^\infty$ converges to a map
$b_g\colon\pG\setminus \{g^-\} \rightarrow G$ uniformly on every
compact subset of $\pG\setminus \{g^-\}$. Moreover, there exists a
constant $C'(g,r)>0$ such that for any $\xi,\eta\in \pG\setminus
V_r(g^-)$ the following estimate holds:
$$\rho(b_g(\xi),b_g(\eta))\leq C'(g,r)\tau^{(\xi|\eta)}.$$
\elem

\bpf Let $\xi,\eta\in \pG\setminus V_r(g^-)$. By the cocycle
identity, we have
$$b_{g,n+1}(\xi)=b_{g,n}(\xi)c(g^{-1},g^n\xi).$$
Thus, thanks to Lemma ~\ref{Livlem1}, we have
\begin{eqnarray*}
\rho(b_{g,n}(\xi),b_{g,n+1}(\xi))&=&\rho(e,c(g^{-1},g^n\xi))\\
&=&\rho(c(g^{-1},g^n  g^+),c(g^{-1},g^n \xi))\\
&\leq& A_{g^{-1}}\tau^{|g^{-n}|+(g^+|\xi)-C(g,r)}\\
&\leq& A_{g^{-1}}\tau^{|g^{-n}|-C(g,r)}.
\end{eqnarray*}
Since $g$ is an infinite order element, there exist constants
$s,t>0$ such that $|g^{-n}|\geq sn-t$ for any $n\in \0 N$ (see
\cite[Chapitre 8,~ Proposition 21]{gh}). Thus the first statement
holds.

Since the above cocycle identity and Lemma ~\ref{Livlem1} imply
$$\rho(b_{g,n}(\xi),b_{g,n}(\eta))\leq \sum_{k=0}^{n-1}
\rho(c(g^{-1},g^k\xi),c(g^{-1},g^k\eta))\leq \sum_{k=1}^{n-1}
A_{g^{-1}}\tau^{|g^{-k}|+(\xi|\eta)-C(g,r)},$$
the second statement holds with
$$C'(g,r)=\sum_{k=0}^\infty A_{g^{-1}}\tau^{|g^{-k}|-C(g,r)}.$$
\epf

\blem \label{Livlem4} Let $g\in \Gamma$ be an infinite order element.
\begin{itemize}
\item [(1)] Let $h\in \Gamma$ with $h g^+\neq g^{-}$. Then $b_g(h
g^+)=c(h^{-1},g^+)^{-1}$ holds. \item [(2)] Let $k\in \Gamma$ and
$\xi\in \pG\setminus \{g^-,k g^-\}$. Then
$b_g(\xi)=b_{kgk^{-1}}(\xi)c(k^{-1},g^+)^{-1}$ holds.
\end{itemize}
\elem

\bpf (1). Lemma ~\ref{Livlem2} and Lemma~\ref{Livlem3} imply
\begin{eqnarray*}
b_g(h g^{+})&=&\lim_{n\to\infty}b_{g,n}(h(g^nh)^+)
=\lim_{n\to \infty}c(g^{-n},h(g^nh)^+)\\
&=&\lim_{n\to \infty}c(h^{-1},(g^nh)^+)^{-1}c(h^{-1}g^{-n},(g^nh)^+)\\
&=&c(h^{-1},g^+)^{-1}.
\end{eqnarray*}

(2). Let $l\in \Gamma$ with $l g^+\in \pG\setminus \{g^-,k g^-\}$.
Then (1) implies
\begin{eqnarray*}b_g(l g^+)&=&c(l^{-1},g^+)^{-1}
=c(kl^{-1},k g^+)^{-1}c(k^{-1},g^+)^{-1}\\
&=&c(kl^{-1},(kgk^{-1})^+)^{-1}c(k^{-1},g^+)^{-1}\\
&=&b_{kgk^{-1}}(lk^{-1}(kgk^{-1})^+)c(k^{-1},g^+)^{-1}\\
&=&b_{kgk^{-1}}(l g^+)c(k^{-1},g^+)^{-1}.
\end{eqnarray*}
Since $b_g$ and $b_{kgk^{-1}}$ are continuous on $\pG\setminus
\{g^-,k g^-\}$ and $\Gamma g^+$ is dense in $\pG$, we get the
statement. \epf

\bpff{Proof of Theorem~\ref{Liv}} We fix an infinite order element
$g\in \Gamma$. Then Lemma~\ref{Livlem4},(2) with $k\in \Gamma$
satisfying $g^-\neq k g^-$ shows that $b_g$ has a unique
continuous extension $b\colon\pG\rightarrow G$, which satisfies
$b(g^-)=b_{kgk^{-1}}(g^-)c(k^{-1},g^+)^{-1}$. Note that this value
does not depend on the choice of $k$ as above.

We first show the H\"older continuity of $b$. We take $h_1,h_2\in
\Gamma$ such that $\omega_0:=g^-$, $\omega_1:=h_1 g^-$ and
$\omega_2:=h_2 g^-$ are distinct points. We choose $r>0$
satisfying $V_r(\omega_i)\cap V_r(\omega_j)=\emptyset$ for any
$i\neq j$ and set
$$A=\max\{C'(g,r),C'(h_1gh_1^{-1},r),C'(h_2gh_2^{-1},r)\}.$$
Note that for any $\xi\neq \eta \in \pG$, there exists
$i\in \{0,1,2\}$ such that $\xi,\eta\in \pG\setminus V_r(\omega_i)$.
Now Lemma~\ref{Livlem3} and Lemma~\ref{Livlem4},(2) imply
$$\rho(b(\xi),b(\eta))\leq A\tau^{(\xi|\eta)}.$$

We claim that $b(h g^+)=c(h^{-1},g^+)^{-1}$ holds for any $h\in
\Gamma$, which has already been shown for $h$ with $h g^+ \neq
g^-$. Assume that $h\in \Gamma$ satisfies $h g^+=g^-$. Then
$g^-\neq h g^-$ and
$$b(h g^+)=b(g^-)=b_{hgh^{-1}}(g^-)c(h^{-1},g^+)^{-1}.$$
Since we have $(hgh^{-1})^+=g^-$, the claim follows from
$$b_{hgh^{-1}}(g^-)=b_{hgh^{-1}}((hgh^{-1})^+)
=\lim_{n\to\infty}c((hgh^{-1})^{-n},(hgh^{-1})^+)=e.$$

For any $h,k\in \Gamma$, we have
$$b(k g^+)b(h^{-1}k g^+)^{-1}=
c(k^{-1},g^+)^{-1}c(k^{-1}h,g^+)=c(h,k g^+).$$ Since $c(h,\cdot)$
and $b$ are continuous and the $\Gamma$-orbit of $g^+$ is dense in
$\pG$, this finishes the proof. \epf
\bigskip
\section{Types of harmonic measures}
Let $g\in \Gamma$ be an infinite order element. We define
$$
r(g)=K(g^{-1},g^+).
$$
Note that being a nonzero positive harmonic function,
$K(\cdot,\xi)$ is nowhere vanishing, so that $r(g)>0$.
We can also write
\begin{eqnarray*}
r(g)&=&\lim_{n\to+\infty}\frac{G(e,g^{n+1})}{G(e,g^n)}
=\lim_{n\to+\infty}G(e,g^n)^{1/n} \\
&=&\lim_{n\to+\infty}\frac{F(e,g^{n+1})}{F(e,g^n)}
=\lim_{n\to+\infty}F(e,g^n)^{1/n}\\
&=&\sup_{n\ge1}F(e,g^n)^{1/n},
\end{eqnarray*}
where the last equality follows from
$F(e,g^n)F(e,g^m)=F(e,g^n)F(g^n,g^{m+n})\le F(e,g^{m+n})$. Since
$F(e,g^n)\le1$, we see that $r(g)\le1$. We put $r(g)=1$ for any
finite order element $g\in\Gamma$.

\blem The function $r$ on $\Gamma$ is a class function satisfying
$r(g^k)=r(g)^k$ for $k\in\mathbb{N}$. If $\mu$ is symmetric, then
$r(g)=r(g^{-1})$. \elem

\bpf The Martin kernel is a cocycle, that is,
$$
K(gh,\xi)=K(g,\xi)K(h,g^{-1}\xi),\ \ \hbox{and thus} \ \
K(h^{-1},\xi)=K(h,h\xi)^{-1}.
$$
Using that $(g^k)^+=g^+$ and $(hgh^{-1})^+=hg^+$ one then easily
checks that $r(g^k)=r(g)^k$ for $k\in\0N$ and $r(hgh^{-1})=r(g)$.

When $\mu$ is symmetric, we have $G(e,g^n)=G(e,g^{-n})$. Hence
$r(g)=r(g^{-1})$. \epf

\bex Consider the simple random walk defined by the canonical
symmetric generating set $S$ of $\mathbb{F}_N$. Then
$$
F(e,s)=\frac{1}{2N-1}
$$
for $s\in S$, see e.g.~\cite[Sect. 2a]{l}. It follows that
$$
F(e,g)=(2N-1)^{-|g|}
$$
for any $g\in\0F_N$. We can then conclude that
$$
r(g)=(2N-1)^{-\ell_g},
$$
where $\ell_g$ is the minimal length of elements in the conjugacy
class of $g$. \eex

\blem
If $g\in\Gamma$ is an infinite order element, then $r(g)<1$.
\elem

\bpf Let $|g|=l$. By \cite[Chapitre 8, Proposition 21]{gh}, there
exists a $(\lambda,c)$-quasi-geodesic ray $\{f(n)\}_{n=0}^\infty$
on $\Gamma$ such that $f(ln)=g^n$ for all $n\geq 0$. By
quasi-geodesic stability \cite[Chapitre~5, Th\'eor\`eme 6]{gh},
there exists a positive constant $H=H(\delta,\lambda,c)$ such that
the quasi-geodesic segment $f([lm,ln])$ is in the
$H$-neighbourhood of some geodesic segment $\alpha_{m,n}\in
[g^m,g^n]$ for $m<n$. Therefore, thanks to Corollary~\ref{MultC2},
there exists $C\geq 1$ such that
$$F(e,g^{m+n})\leq CF(e,g^m)F(g^m,g^{m+n})$$
for $m,n\in\mathbb{N}$. This implies
$$F(e,g^{m(n+1)})\leq CF(e,g^{mn})F(g^{mn},g^{m(n+1)})
=CF(e,g^{mn})F(e,g^m).$$
Hence we have
$$r(g^m)=\lim_{n\to\infty}\frac{F(e,g^{m(n+1)})}{F(e,g^{mn})}
\leq CF(e,g^m).$$ Since $F(e,g^m)\to0$ as $m\to\infty$ by
(\ref{1eExp}), we see that $r(g^m)<1$ for sufficiently large $m$.
As $r(g^m)=r(g)^m$ by the previous lemma, it follows that
$r(g)<1$.\epf

Denote by ${\mcal R}(\Gamma,\mu)={\mcal
R}(\Gamma,\partial\Gamma,\nu)$ the orbit equivalence relation
defined by the action of $\Gamma$ on $(\partial\Gamma,\nu)$.

Recall that by definition the ratio set
$r(\Gamma,\partial\Gamma,\nu)$ consists of all $\lambda\ge0$ such
that for any $\eps>0$ and any subset $A\subset\partial\Gamma$ of
positive measure there exists $g\in\Gamma$ such that the set of
$\omega\in gA\cap A$ satisfying
$$
\left|\frac{dg\nu}{d\nu}(\omega)-\lambda\right|<\eps
$$
has positive measure.

Note that ${\mcal R}(\Gamma,\mu)$ is ergodic, amenable and of type
III (by~\cite[Theorem~3.2.1]{k1}). Hence $\{0,1\}\subset
r(\Gamma,\partial\Gamma,\nu)$, and
$r(\Gamma,\partial\Gamma,\nu)\setminus\{0\}$ is a closed
multiplicative subgroup of $(0,+\infty)$. One says that ${\mcal
R}(\Gamma,\mu)$ is of type III$_0$, III$_\lambda$ ($0<\lambda<1$)
or III$_1$ depending on whether this group is $\{1\}$,
$\{\lambda^n\}_{n\in\0Z}$ or $(0,+\infty)$. Recall also that for
$0<\lambda\le1$ there is only one amenable ergodic equivalence
relation of type~III$_\lambda$.

We can now formulate our main result.

\bthm \label{Main} Let $\Gamma$ be a nonelementary hyperbolic
group, $\nu$ the harmonic measure on $\partial\Gamma$ defined by a
finitely supported nondegenerate probability measure $\mu$ on
$\Gamma$. Then $r(\Gamma,\partial\Gamma,\nu)\setminus\{0\}$ is the closed
multiplicative subgroup of $(0,+\infty)$ generated by $\{r(g)\}_{g\in \Gamma}$.
In particular, ${\mcal R}(\Gamma,\mu)$ is never of type~III$_0$. \ethm

Note that since the harmonic measure is nonatomic and any infinite
order element has only two fixed points, the crossed product
$L^\infty(\partial\Gamma,\nu)\rtimes\Gamma$ is a factor if
$\Gamma$ is torsion-free. The theorem gives then the type of this
factor.

To prove the theorem we need the following lemmas.

\blem \label{FreeL} Let $g\in\Gamma$ be an infinite order element.
Then there exists $L_1>0$ such that for any $x\in\Gamma$, there is
$y\in\Gamma$ such that $|y|\leq L_1$ and $|xyg^n|\geq
|x|+|g^n|-L_1$ for all $n\in\mathbb{Z}$. \elem

\bpf Since $\Gamma$ is nonelementary, there is $a\in\Gamma$ such
that for $h=aga^{-1}$ we have
$\{g^{\pm}\}\cap\{h^{\pm}\}=\emptyset$, see e.g. the proof
of~\cite[Chapitre 8, Th\'eor\`eme 37]{gh}. Then there is $M>0$
such that
$$
(g^n|h^m)\leq M\ \ \hbox{for all}\ n,m\in\mathbb{Z},
$$
since otherwise we could find a subsequence $\{x_n\}^\infty_{n=1}$
of $\{g^n\}_{n\in\0Z}$ converging to a point in $\{g^\pm\}$ and a
subsequence $\{y_n\}^\infty_{n=1}$ of $\{h^n\}_{n\in\0Z}$
converging to a point in $\{h^\pm\}$ such that $\{x_n\}_n$ and
$\{y_n\}_n$ are equivalent.

Take $L_1\geq2M+2\delta+3|a|$. Let $x\in \Gamma$. If
$(x^{-1}|g^n)\leq M+\delta$ for any $n\in\mathbb{Z}$, then we have
$$2M+2\delta\geq |x|+|g^n|-|xg^n|$$
and hence we can take $y=e$. If there is
$m\in\mathbb{Z}$ such that $(x^{-1}|g^m)>M+\delta$, then for any
$n\in\mathbb{Z}$,
$$M\geq(h^n|g^m)\geq\min\{(x^{-1}|h^n),
(x^{-1}|g^m)\}-\delta=(x^{-1}|h^n)-\delta.$$
Hence we obtain
$$M+\delta\geq(x^{-1}|h^n)=(x^{-1}|ag^na^{-1}).$$
Therefore we have
$$2M+2\delta\geq|x|+|ag^na^{-1}|-|xag^na^{-1}|
\geq|x|+|g^n|-|xag^n|-3|a|.$$
Thus in this case we can take $y=a$.
\epf

\blem \label{Lambda} If $\{r(g)\}_{g\in \Gamma}$ is a subset of
$\{\lambda^n\}_{n\in \0Z}$ for some $0<\lambda<1$, then
$r(\Gamma,\partial\Gamma,\nu)\setminus\{0\}$ is a subgroup of $
\{\lambda^n\}_{n\in \0Z}$. \elem

\bpf Let $T=-2\pi/\log \lambda$ and set
$c(g,\xi)=K(g,\xi)^{\sqrt{-1}T}$. Then thanks to
Corollary~\ref{HolderC}, the cocycle $c$ satisfies the assumption
of Theorem~\ref{Liv} with $G=\0 T=\0 R/2\pi \0 Z$ and there exists
a continuous map $b\colon\pG\rightarrow \0 T$ satisfying
$c(g,\xi)=b(\xi)b(g^{-1}\xi)^{-1}$ for all $g\in \Gamma$ and
$\xi\in \pG$. We choose a Borel map $f\colon\pG\rightarrow
(\lambda,1]$ satisfying $b(\xi)=f(\xi)^{\sqrt{-1}T}$ for all $\xi$
and set $\nu_1=f\nu$. Then $\nu_1$ is equivalent to $\nu$ and
satisfies
$$\frac{dg\nu_1}{d\nu_1}(\omega)\in \{\lambda^n\}_{n\in \0Z},\quad
\forall g\in \Gamma,\; \forall \omega\in \pG.$$
This shows the statement.
\epf

\bpff{Proof of Theorem~\ref{Main}} The proof is inspired by
Bowen's computation of the ratio set of a Gibbs measure
in~\cite[Lemma 8]{b}.

Thanks to Lemma~\ref{Lambda}, to prove the theorem it suffices to
show that $r(g)$ belongs to the ratio set for any $g\in\Gamma$. So
let $g\in\Gamma$ be an infinite order element and
$\lambda=r(g)=K(g^{-1},g^+)$. Let $\varepsilon>0$ and
$A\subset\partial\Gamma$ with $\nu(A)>0$.

There exists a visual metric $\rho$ on $\partial\Gamma$ such that
for almost every $\omega\in A$ we have
$$
\lim_{r\to0}\frac{\nu(A\cap B(\omega,r))}{\nu(B(\omega,r))}=1.
$$
Indeed, by~\cite[Theorems~9.1 and~9.2]{bs} there exists a visual
metric $\rho$ such that $(\partial\Gamma,\rho)$ embeds
isometrically into $\0R^n$ for some $n$. Then the above
convergence holds by a classical result of Besicovitch, see
e.g.~\cite{f}. Note that using Theorem~\ref{Gibbs} one can
then show that the convergence holds for any visual metric,
but we do not need this.

Therefore there exist $\omega\in A$, $N_0\in\0N$ and closed
neighbourhoods $W(R)$ of $\omega$ for $R\ge N_0$ such that
\begin{equation} \label{eMain0}
U(\omega,R-N_0)\supset W(R)\supset U(\omega,R)
\end{equation}
and
\begin{equation} \label{eMain1}
\frac{\nu(A\cap W(R))}{\nu(W(R))}\to 1 \ \ \hbox{as}\ \
R\to\infty.
\end{equation}
Fix a geodesic ray $\{w(n)\}_{n=0}^\infty\in[e, \omega]$.

By Corollary~\ref{HolderC} there is $L_0>0$ such that
\begin{equation} \label{eMain2}
e^{-\varepsilon}\leq\frac{K(h,\xi)}{K(h,\eta)} \leq
e^{\varepsilon} \ \ \mbox{if}\ \ (\xi|\eta)-|h|\geq L_0.
\end{equation}

By Lemma~\ref{FreeL} there exists $L_1>0$ such that for any $x$
there is $y\in\Gamma$ such that
\begin{equation} \label{eMain3}
|y|\leq L_1\ \ \hbox{and}\ \ |xyg^n|\geq |x|+|g^n|-L_1\ \
\mbox{for}\ \ n\in\mathbb{Z}.
\end{equation}

Since $K(g^{-1},\cdot)$ is continuous, there is $L_2>0$ such that
\begin{equation} \label{eMain4}
\lambda e^{-\varepsilon}\leq K(g^{-1},\zeta) \leq \lambda
e^{\varepsilon}\ \ \mbox{for}\ \ \zeta\in U(g^+,L_2).
\end{equation}

Choose an integer $L>\max\{\delta,|g|,L_0,L_1,L_2\}$.

\smallskip

Fix $N\in\0N$ to be specified later and  put $x=w(N)$. Then
choose $y$ satisfying~(\ref{eMain3}) and put $h=xy$. Consider
$V=U(hg^+,N+14L)$. We claim that
\begin{equation} \label{eMain5}
h^{-1}V\subset U(g^+,11L)\ \ \mbox{and}\ \ gh^{-1}V\subset
U(g^+,10L).
\end{equation}
Indeed, let $\zeta\in U(hg^+,N+14L)$. Then $(\zeta|hg^+)>N+14L$.
Since $|h|<N+L$, we thus get
$$
(h^{-1}\zeta|g^+)\ge (\zeta|hg^+)-|h|>13L.
$$
Hence $h^{-1}\zeta\in U(g^+,13L-2\delta)\subset U(g^+,11L)$. Since
$|gh^{-1}|\leq N+2L$, we similarly get $gh^{-1}V\subset
U(g^+,10L)$.

Next we claim that
\begin{equation} \label{eMain6}
hU(g^+,10L)\subset U(hg^+,N+4L).
\end{equation}
Let $\zeta\in U(g^+,10L)$. Take geodesic rays $z\in[e,\zeta]$ and
$v\in[e,g^+]$. Note that
$$\liminf_{n\to\infty}(z(n)|g^n)
\geq\lim_{n\to\infty}(z(n)|v(n))-\delta>9L.$$ Since by
(\ref{eMain3})
$$
(h^{-1}|g^n)=\frac12(|xy|+|g^n|-|xyg^n|)\le\frac12(|xy|+L-|x|)\leq L
$$
for any $n\in\mathbb{Z}$, we get
$$L\geq (h^{-1}|g^n)\geq\min\{(h^{-1}|z(n)), (z(n)|g^n)\}-\delta,$$
and so $(h^{-1}|z(n))\leq L+\delta\leq 2L$ for $n$ large enough,
so that, since $|h|\ge N-L$,
$$|hz(n)|\geq |h|+|z(n)|-4L\geq |z(n)|+N-5L.$$
Therefore using that $|hg^n|\ge|g^n|+N-L$ by (\ref{eMain3}), we
get
\begin{eqnarray*}
(hz(n)|hg^n)
&\geq&\frac{1}{2}(|z(n)|+N-5L+|g^n|+N-L-d(z(n),g^n))\\
&=&N-3L+(z(n)|g^n) \\
&>&N+6L
\end{eqnarray*}
for sufficiently large $n$. Hence $(h\zeta|hg^+)>N+6L$ and thus
$h\zeta\in U(hg^+,N+6L-2\delta)\subset U(hg^+,N+4L)$.

By (\ref{eMain5}) and (\ref{eMain6}), we have
\begin{equation} \label{eMain7}
hgh^{-1}V\subset U(hg^+,N+4L).
\end{equation}

For any $\zeta\in V$, by using the cocycle property we obtain
$$K(hg^{-1}h^{-1},\zeta)
=K(g^{-1},h^{-1}\zeta)\frac{K(h,\zeta)}{K(h,hgh^{-1}\zeta)}.$$
Since $h^{-1}\zeta\in U(g^+,11L)$ by~(\ref{eMain5}), the first
factor on the right hand side of the above equality is in
$[\lambda e^{-\varepsilon},\lambda e^{\varepsilon}]$
by~(\ref{eMain4}). Since $|h|\leq N+L$, we have,
by~(\ref{eMain7}),
$$(\zeta|hgh^{-1}\zeta)
\geq\min\{(\zeta|hg^+),(hgh^{-1}\zeta|hg^+)\}-2\delta>N+4L-2\delta
\geq N+2L \geq |h|+L.$$ Hence the second factor is in
$[e^{-\varepsilon},e^{\varepsilon}]$ by~(\ref{eMain2}). Thus
\begin{equation} \label{eMain8}
\lambda e^{-2\varepsilon}\leq K(hg^{-1}h^{-1},\zeta)
\leq \lambda e^{2\varepsilon}\ \ \mbox{for}\ \ \zeta\in V.
\end{equation}
To complete the proof it suffices to show that by choosing
sufficiently large $N$ we could arrange $\nu(V\cap A\cap
hg^{-1}h^{-1}A)>0$.

\smallskip

We shall check first that
\begin{equation} \label{eMain9}
V\subset U(\omega,N-4L)\ \ \mbox{and}\ \ hgh^{-1}V\subset
U(\omega,N-4L).
\end{equation}
By virtue of (\ref{eMain7}) it suffices to show that
$U(hg^+,N+4L)\subset U(\omega,N-4L)$. Since
\begin{eqnarray*}
2(w(n)|hg^n)&\ge&n+|hg^n|-|w(n)^{-1}w(N)|-|w(N)^{-1}hg^n| \\
&\geq&n+(N+|g^n|-L)-(n-N)-(L+|g^n|) \\
&=&2N-2L,
\end{eqnarray*}
we have $(\omega|hg^+)\ge N-L$, so that $\omega\in
U(hg^+,N-3L)\subset U(hg^+,N-4L+\delta)$. By~(\ref{eGibbs}), it
follows that $U(hg^+,N-3L)\subset U(hg^+,N-4L+\delta)\subset
U(\omega,N-4L)$, and (\ref{eMain9}) is proved. Similarly we have
$U(\omega,N-3L)\subset U(hg^+,N-4L)$.

Then by Theorem~\ref{Gibbs}, for $N>4L$ and
$\{v(n)\}^\infty_{n=0}\in[e,hg^+]$ we obtain
\begin{eqnarray*}
\nu(V)&\geq&\frac{F(e,v(N+14L))}{D} \\
&\geq&\frac{F(e,v(N-4L))F(v(N-4L),v(N+14L))}{D} \\
&=&\frac{F(e,v(N-4L))F(e,v(N-4L)^{-1}v(N+14L))}{D} \\
&\geq&\frac{\nu(U(hg^+,N-4L))F(e,v(N-4L)^{-1}v(N+14L))}{D^2} \\
&\geq&\frac{\nu(U(\omega,N-3L))F(e,v(N-4L)^{-1}v(N+14L))}{D^2}
\end{eqnarray*}
and similarly
$$
\nu(U(\omega,N-3L)) \ge\frac{\nu(U(\omega,N-N_0-4L))
F(e,w(N-N_0-4L)^{-1}w(N-3L))}{D^2}.
$$
Thus if we put
$$c=\min\left\{\left.\frac{F(e,x_1)F(e,x_2)}{D^4}\right|
|x_1|=18L, |x_2|=N_0+L\right\}>0,$$ then
$$\nu(V)\geq c\nu(U(\omega,N-N_0-4L))\ge c\nu(W(N-4L)),$$
where the second inequality follows from~(\ref{eMain0}). Let
$\varepsilon_N\ge0$ be such that
$$\nu(A\cap W(N-4L))=(1-\varepsilon_N)\nu(W(N-4L)).$$
Then $\eps_N\to0$ as $N\to\infty$ by (\ref{eMain1}). So if we
denote $A\cap V$ by $X$, then since $V\subset W(N-4L)$
by~(\ref{eMain0}) and~(\ref{eMain9}), we get
\begin{eqnarray*}
\nu(X)
&\geq&\nu(A\cap W(N-4L))-\nu(W(N-4L)\setminus V) \\
&=&(1-\varepsilon_N)\nu(W(N-4L))
-(\nu(W(N-4L))-\nu(V)) \\
&=&-\varepsilon_N\nu(W(N-4L))+\nu(V) \\
&\geq&(c-\varepsilon_N)\nu(W(N-4L))>0
\end{eqnarray*}
for $N$ sufficiently large. Since
$$\nu(hgh^{-1}X)=\int_XK(hg^{-1}h^{-1},\zeta)d\nu(\zeta)
\geq\lambda e^{-2\varepsilon}\nu(X) \geq(c-\varepsilon_N)\lambda
e^{-2\varepsilon}\nu(W(N-4L)),$$ we similarly obtain
$$\nu((hgh^{-1}X)\cap A)
\geq\{(c-\varepsilon_N)\lambda e^{-2\varepsilon}-\varepsilon_N\}
\nu(W(N-4L)).$$ The latter expression is strictly positive if $N$
is sufficiently large. Thus $\nu(V\cap A\cap
hg^{-1}h^{-1}A)=\nu(X\cap hg^{-1}h^{-1}A)>0$, and the proof of the
theorem is complete. \epf

\bigskip

\flushleft{Masaki Izumi, Department of Mathematics, Graduate
School of Science, Kyoto University, Sakyo-ku, Kyoto 606-8502,
Japan\\
{\it e-mail}: izumi@math.kyoto-u.ac.jp}

\flushleft{Sergey Neshveyev, Mathematics Institute, University of
Oslo, PB 1053 Blindern, Oslo 0316, Norway\\
{\it e-mail}: sergeyn@math.uio.no}

\flushleft{Rui Okayasu, Department of Mathematics Education, Osaka
Kyoiku University, Kashiwara, Osaka 582-8582, Japan\\
{\it e-mail}: rui@cc.osaka-kyoiku.ac.jp}

\end{document}